\documentclass[conference]{IEEEtran}
\usepackage[cmex10]{amsmath}
\usepackage{graphicx}
\usepackage{acronym}
\usepackage{amsmath}
\usepackage{amsfonts}
\usepackage{xspace}
\usepackage{pgfplots}
\usepackage[noend]{algpseudocode}
\usepackage{algorithm}
\usepackage{algorithmicx}

\algrenewcommand\algorithmicrequire{\textbf{Input:}}
\algrenewcommand\algorithmicensure{\textbf{Input:}}

\newcommand{\dt}{\Delta t}

\title{\Large \bf Optimization Outperforms Unscented Techniques for Nonlinear Smoothing} 
\author{\IEEEauthorblockN{
Payton Howell \quad \quad
Aleksandr Aravkin 
}\\
\IEEEauthorblockN{Applied Mathematics, University of Washington, Seattle, WA USA}
}

\usepackage[pscoord]{eso-pic}
\newcommand{\placetextbox}[3]{
 \setbox0=\hbox{#3}
 \AddToShipoutPictureFG*{ \put(\LenToUnit{#1\paperwidth},\LenToUnit{#2\paperheight}){\vtop{{\null}\makebox[0pt][c]{#3}}}
 }
 }

\pgfplotsset{compat=1.18}
\begin{document}

 \placetextbox{.23}{0.055}{l}{\small{\copyright\,\copyright 2025 IEEE. Personal use of this material is permitted. Permission from IEEE must be obtained for all other uses, in any current or future media, including reprinting/republishing this material for advertising or promotional purposes, \\ creating new collective works, for resale or redistribution to servers or lists, or reuse of any copyrighted component of this work in other works.}}

\maketitle

% use for special paper notices
%\IEEEspecialpapernotice{(Invited Paper)}

% make the title area
%\maketitle

\begin{abstract}
We review optimization-based approaches to smoothing nonlinear dynamical systems. 
These approaches leverage the fact that the Extended Kalman Filter and corresponding 
smoother can be framed as the Gauss-Newton method for a nonlinear least squares maximum a posteriori loss, and stabilized with standard globalization techniques. We compare the performance of the Optimized Kalman Smoother (OKS) to Unscented Kalman smoothing  techniques, and show that they achieve significant improvement for highly nonlinear systems, particularly in noisy settings. The comparison is performed across standard parameter choices (such as the trade-off between process and measurement terms). To our knowledge, this is the first comparison of these methods in the literature. 
\end{abstract}

%\vspace{-0.2in}

% my_acronyms.sty
%
% define acronyms used with my thesis

% #'s
\acrodef{2D}{two-dimensional}
\acrodef{3D}{three-dimensional}

% A
\acrodef{ABE}{Autonomous Benthic Explorer}
\acrodef{ADCP}{acoustic Doppler current profiler}
\acrodef{ADV}{acoustic Doppler velocimeter}
\acrodef{Alvin}{Alvin}
\acrodef{AFRL}{Air Force Research Laboratory}
\acrodef{AHRS}{attitude and heading reference system}
\acrodef{Autosub}{Autosub}
\acrodef{AUG}{autonomous underwater glider}
\acrodef{AUV}{autonomous underwater vehicle}
\acrodef{AGU}{American Geophysical Union}
\acrodef{AON}{Arctic Observing Network}
\acrodef{AOPE}{Applied Ocean Physics and Engineering}
\acrodef{AS}{asymptotically stable}
\acrodef{ACFR}{Australian Centre for Field Robotics}
\acrodef{ASV}{autonomous surface vehicle}

% B
\acrodef{BCC}{brightness constancy constraint}
\acrodef{BIO}{Bedford Institute of Oceanography}

% C
\acrodef{CAD}{computer aided design}
\acrodef{CenSSIS}{Center for Subsurface Sensing and Imaging Systems}
\acrodef{CCD}{charge coupled device}
\acrodef{CI}{covariance intersection}
\acrodef{CG}{conjugate gradients}
\acrodef{COG}{course over ground}
\acrodef{CPU}{central processing unit}
\acrodef{CT}{continuous-time}
\acrodef{CB}{center of buoyancy}
\acrodef{CG}{center of gravity}
\acrodef{CGSN}{Canadian Gravity Standardization Net}
\acrodef{COM}{center of mass}
\acrodef{CO2}{carbon dioxide}
\acrodef{CPS}{control power supply}
\acrodef{CSAC}{chip-scale atomic clock}
\acrodef{CST}{{\it IEEE Transactions on Control Systems Technology}}
\acrodef{CSV}{Comma Separated Values}
\acrodef{CTD}{conductivity-temperature-depth}
\acrodef{CISE}{Computer and Information Science and Engineering}

% D
\acrodef{DOF}{degree of freedom}
\acrodef{DoD}{Department of Defense}
\acrodef{DSL}{Deep Submergence Laboratory}
\acrodef{DT}{discrete-time}
\acrodef{DVL}{Doppler velocity log}
\acrodef{DPM}{digital panel meter}
\acrodef{DR}{dead reckoning}
\acrodef{DRI}{department research initiative}
\acrodef{DWH}{Deepwater Horizon}

% E
\acrodef{ESDF}{exactly sparse delayed-state filter}
\acrodef{EM}{electromagnetic}
\acrodef{EKF}{extended Kalman filter}
\acrodef{EIF}{extended information filter}
\acrodef{ERC}{Engineering Research Center}
\acrodef{EPR}{East Pacific Rise}
\acrodef{EPSL}{{\it Earth and Planetary Science Letters}}

% F
\acrodef{FastSLAM}{Factored Solution to SLAM}
\acrodef{FOG}{fiber optic gyro}
\acrodef{FOV}{field of view}
\acrodef{FFT}{fast Fourier transform}
\acrodef{FBD}{free body diagram}
\acrodef{FIRST}{For Inspiration and Recognition of Science and Technology}

% G
\acrodef{GMRF}{Gaussian Markov random field}
\acrodef{GPS}{global positioning system}
\acrodef{GAS}{globally asymptotically stable}
\acrodef{GA}{geometric algebra}
\acrodef{GUI}{graphical user interface}

% H
\acrodef{HUGIN}{HUGIN}
\acrodef{HMMV}{H\r{a}kon Mosby Mud Volcano}
\acrodef{HOV}{human occupied vehicle}
\acrodef{HROV}{hybrid remotely operated vehicle}
\acrodef{HTF}{Hydrodynamics Test Facility}

% I
\acrodef{ICRA}{{\it IEEE International Conference on Robotics and Automation}}
\acrodef{IF}{information filter}
\acrodef{IFREMER}{French Institute for the Research and Exploitation of the Sea}
\acrodef{IFE}{Institute for Exploration}
\acrodef{INU}{inertial navigation unit}
\acrodef{INS}{inertial navigation system}
\acrodef{IR}{infrared}
\acrodef{IMU}{inertial measurement unit}
\acrodef{INS}{inertial navigation system}
\acrodef{IBCAO}{International Bathymetric Chart of the Arctic Ocean}
\acrodef{IROS}{{\it IEEE International Conference on Robotics and Intelligent Systems}}
\acrodef{iUSBL}{inverted USBL}
\acrodef{OWTTiUSBL}{One-Way Travel-Time inverted USBL}

% J
\acrodef{Jason}{Jason}
\acrodef{JFR}{{\it Journal of Field Robotics}}
\acrodef{JdF}{Juan de Fuca}
\acrodef{JHU}{Johns Hopkins University}
\acrodef{JPEG}{Joint Photographic Experts Group}

% K
\acrodef{KF}{Kalman filter}
\acrodef{KYP}{Kalman-Yakubovich-Popov}
\acrodef{KISS}{Keck Institute for Space Studies}

% L
\acrodef{LADCP}{lowered acoustic Doppler current profiler}
\acrodef{LBL}{long-baseline}
\acrodef{LCM}{lightweight communication and marshaling}
\acrodef{LG}{linear Gaussian}
\acrodef{LKY}{Lefschetz-Kalman-Yakubovich}
\acrodef{LMedS}{least median of squares}
\acrodef{LLE}{Linear Lyapunov Equation}
\acrodef{LQR}{Linear Quadratic Regulation}
\acrodef{LQGR}{Linear Quadratic Gaussian Regulation}
\acrodef{LSSL}{Louis S. St Laurent}
\acrodef{LRAUV}{long-range AUV}

% M
\acrodef{MAP}{maximum \emph{a posteriori}}
\acrodef{MBARI}{Monterey Bay Aquarium Research Institute}
\acrodef{MBN}{mosaic-based navigation}
\acrodef{MEF}{Main Endeavor Field}
\acrodef{MIT}{Massachusetts Institute of Technology}
\acrodef{MLE}{maximum likelihood estimate}
\acrodef{MRF}{Markov random field}
\acrodef{MIZ}{Marginal Ice Zone}
\acrodef{MATE}{Marine Advanced Technology Education}
\acrodef{MCR}{Mid-Cayman Rise}
\acrodef{MEMS}{micro-electrical-mechanical systems}
\acrodef{MBARI}{Monterey Bay Aquarium Research Institute}
\acrodef{MIMO}{multiple-input, multiple-output}
\acrodef{MOR}{mid-ocean ridge}
\acrodef{MVCO}{Martha's Vineyard Coastal Observatory}
\acrodef{MCL}{mission critical level}
\acrodef{MO}{mission objective}
\acrodef{MSD}{mass spring damper}
% N
\acrodef{NavEst}{navigation estimation}
\acrodef{NDSEG}{National Defense Science and Engineering Graduate}
\acrodef{NEES}{normalized estimation error squared}
\acrodef{NDSF}{National Deep Submergence Facility}
\acrodef{NMEA}{National Marine Electronics Association}
\acrodef{NSF}{National Science Foundation}
\acrodef{NTP}{network time protocol}
\acrodef{NAS}{National Academy of Sciences}
\acrodef{NDSF}{National Deep Submergence Facility}
\acrodef{NHS}{Natick High School}
\acrodef{NLO}{nonlinear observer}
\acrodef{NSTA}{National Science Teachers Association}
\acrodef{NSIDC}{National Snow and Ice Data Center}
\acrodef{NUI}{Nereid under-ice}
% O
\acrodef{OS}{operating system}
\acrodef{OWTT}{one-way travel time}
\acrodef{ODE}{ordinary differential equation}
\acrodef{OL}{Second-Order, Open-Loop Observer}
\acrodef{ONR}{Office of Naval Research}
\acrodef{OOI}{Ocean Observing Initiative}
\acrodef{OWTT}{one-way travel time}

% P
\acrodef{PC}{personal computer}
\acrodef{PCB}{printed circuit board}
\acrodef{PDE}{partial differential equation}
\acrodef{PPS}{pulse per second}
\acrodef{PEG}{parameter error gain}
\acrodef{PHF}{Piccard Hydrothermal Field}
\acrodef{PNAS}{{\it Proceedings of the National Academy of Sciences}}
\acrodef{ppb}{part-per-billion}
\acrodef{ppm}{part-per-million}
\acrodef{PR}{Positive Real}
\acrodef{PROV}{polar remotely operated vehicle}
\acrodef{PI}{principal investigator}

% Q

% R
\acrodef{RAM}{random access memory}
\acrodef{RANSAC}{random sample consensus}
\acrodef{RDI}{RD Instruments}
\acrodef{REMUS}{Remote Environmental Monitoring Unit}
\acrodef{RLG}{ring laser gyroscope}
\acrodef{RF}{radio frequency}
\acrodef{RMS}{Royal Mail Steamship}
\acrodef{ROM}{range only measurement}
\acrodef{ROV}{remotely operated vehicle}
\acrodef{RTC}{real-time clock}
\acrodef{ROV}{remotely operated vehicle}
\acrodef{ROVER}{remotely operated vehicle environmental research}
\acrodef{RNS}{Regional Scale Node}
\acrodef{RSS}{{\it Robotics: Science and Systems}}
\acrodef{RTS}{Rauch-Tung-Striebel}

% S
\acrodef{SeaBED}{SeaBED}
\acrodef{SEIF}{sparse extended information filter}
\acrodef{SIFT}{scale invariant feature transform}
\acrodef{SLAM}{simultaneous localization and mapping}
\acrodef{SNAME}{The Society of Naval Architects and Marine Engineers}
\acrodef{SSD}{sum of squared differences}
\acrodef{SFM}{structure-from-motion}
\acrodef{SO}{special orthogonal}
\acrodef{SOA}{state of the art}
\acrodef{SPR}{Strictly Positive Real}
\acrodef{STEM}{Science, Technology, Engineering, and Mathematics }
\acrodef{SSF}{Summer Student Fellow}
\acrodef{SVO}{Scaler, Velocity Observer}
\acrodef{SVD}{singular value decomposition}
\acrodef{SVP}{sound velocity profile}
\acrodef{SPURS}{Salinity Processes in the Upper Ocean Regional Study}
% T
\acrodef{TDMA}{time division multiple access}
\acrodef{TRL}{technology readiness level}
\acrodef{TRO}{{\it IEEE Transactions on Robotics}}
\acrodef{TJTF}{thin junction-tree filter}
\acrodef{TTL}{transistor-transistor logic}
\acrodef{TXCO}{temperature compensated crystal oscillator}
\acrodef{3D}{three dimensional}
\acrodef{TWTT}{two-way travel time}

% U
\acrodef{USBL}{ultra-short-baseline}
\acrodef{UUV}{unmanned underwater vehicle}
\acrodef{UTC}{Coordinate Universal Time}
\acrodef{UDP}{User Datagram Protocol}
\acrodef{UKF}{Unscented Kalman Filter}
\acrodef{UV}{Underwater Vehicle}
\acrodef{UNOLS}{University National Oceanographic Laboratory System}
\acrodef{USGS}{United States Geological Survey}
\acrodef{USNA}{United States Naval Academy}
\acrodef{UNCLOS}{United Nations Convention on the Law of the Sea}
% V
\acrodef{VAN}{visually augmented navigation}
\acrodef{VNL}{vision numerical library}
\acrodef{VTK}{The Visualization Toolkit}
\acrodef{VIGA}{vehicle induced gravimeter acceleration}

% W
\acrodef{WHOI}{Woods Hole Oceanographic Institution}

% X

% Y
\acrodef{YIP}{Young Investigator Program}

% Z
%NEW COMMAND LIST
% this is a command for superscript and subscript on the left
%\newcommand{\preind}[3]{\;{{\small{#1}}\atop{\small{#2}}}#3}
\newcommand{\preind}[3]{\;{{\tiny{#1}}\atop{\tiny{#2}}}\hspace{-0.05in}#3}
\newcommand{\laplace}{\mathcal{L}}
\newcommand{\rb}[1]{\raisebox{-1.5ex}[0pt]{#1}}

\newcommand{\sentry}[0]{{\it Sentry }}
\newcommand{\nereus}[0]{{\it Nereus }}
\newcommand{\jason}[0]{{\it Jason }}
\newcommand{\alvin}[0]{{\it Alvin }}
\newcommand{\iver}[0]{{\it Iver2 }}
\newcommand{\tioga}[0]{{\it Tioga }}

\newcommand{\order}[1]{\ensuremath{\mathcal{O}(10^{#1})}}

% Eriksen uses the name ``Deepglider.''
\newcommand{\deepGlider}{Deepglider}

% These names suck.
\newcommand{\standardTraj}{conventional\xspace}
\newcommand{\USBLTraj}{deep-profiling\xspace}

\newcommand{\effFactor}{R}
\newcommand{\effFactorDive}{R_{dive}}

\newcommand{\TDive}{T_{dive}}
\newcommand{\Tconv}{T_{c}} %jck -- time a conventional glider spends submerged
\newcommand{\Tdeep}{T_{d}} %jck -- time a deep profiling glider spends submerged

\newcommand{\profileDepth}{z}
\newcommand{\profileHeight}{\Delta_z}
\newcommand{\diveSpeed}{U}
\newcommand{\vertSpeed}{W}
\newcommand{\diveAngle}{\theta}
\newcommand{\pumpEnergy}{E_{pump}}
\newcommand{\batteryEnergy}{B}
\newcommand{\pHotel}{P_{hotel}}
\newcommand{\pNav}{P_{nav}}
\newcommand{\diveEnergyConventional}{E_{dive,o}}
\newcommand{\diveEnergyUSBL}{E_{dive}}
\newcommand{\xenduranceConventional}{T_o}
\newcommand{\xenduranceUSBL}{T}

%commands for the 2015 RI:small proposal -----
%2015/07/27 21:06:10  changing some notation for OWTT-USBL.

%super- and subscripts to frames and alignments.
\newcommand{\vehicle}[0]{v}
\newcommand{\world}[0]{w}
\newcommand{\usbl}[0]{u}
\newcommand{\locallevel}[0]{n}  % local north-east-down.

% generic position
\newcommand{\pos}[2]{\preind{#1}{}{\mathbf p}_{#2}}

% generic rotation
\newcommand{\rot}[2]{\preind{#1}{#2}\mathbf{R}}

%owtt-usbl
\newcommand{\upu}[0]{\pos{\usbl}{\mathrm{ASV}}}
\newcommand{\wpv}[0]{\pos{\locallevel}{\mathrm{ASV}}}
\newcommand{\vpw}[0]{\pos{\vehicle}{\world}}
\newcommand{\vuR}[0]{\rot{\vehicle}{\usbl}}
\newcommand{\uvR}[0]{\rot{\usbl}{\vehicle}}
\newcommand{\wvR}[0]{\rot{\locallevel}{\vehicle}}
\newcommand{\vwR}[0]{\rot{\vehicle}{\locallevel}}
\newcommand{\az}[0]{\alpha}
\newcommand{\el}[0]{\gamma}
\newcommand{\rng}[0]{\Gamma}

%scalars
\newcommand{\x}[0]{x(t)}
\newcommand{\xdot}[0]{\dot{x}(t)}
\newcommand{\xddot}[0]{\ddot{x}(t)}
\newcommand{\xhat}[0]{\hat{x}(t)}
\newcommand{\xhatdot}[0]{\dot{\hat{x}}(t)}
\newcommand{\xhatddot}[0]{\ddot{\xhat}(t)}
\newcommand{\deltax}[0]{\Delta x(t)}
\newcommand{\deltaxdot}[0]{\Delta \dot{x}(t)}
\newcommand{\deltaxddot}[0]{\Delta \ddot{x}(t)}
\newcommand{\none}[0]{l_1 s_1 + l_2 s_2}
\newcommand{\ntwo}[0]{l_3 s_1 + l_4 s_2 + r}

\newcommand{\vel}[0]{v(t)}
\newcommand{\veldot}[0]{\dot{v}(t)}
\newcommand{\velhat}[0]{\hat{v}(t)}
\newcommand{\velhatdot}[0]{\dot{\hat{v}}(t)}
\newcommand{\veldelta}[0]{\Delta \vel}
\newcommand{\veldeltadot}[0]{\dot{\Delta \vel}}
\newcommand{\veldeltaq}[0]{\Delta v_q(t)}

\newcommand{\out}[0]{w(t)}
\newcommand{\outhat}[0]{\hat{w}(t)}
\newcommand{\outdelta}[0]{\Delta w(t)}

\newcommand{\invec}[0]{\bm{u}(t)}

\newcommand{\measnoise}[0]{n(t)}

%uuv plant parameters with NO indices subscripts
\newcommand{\thrust}[0]{\tau(t)}
\newcommand{\mass}[0]{m}
\newcommand{\bouyancy}[0]{b}
\newcommand{\quaddrag}[0]{d_{Q}}
\newcommand{\lindrag}[0]{d_{L}}

%uuv plant parameters with indices subscripts
\newcommand{\veli}[0]{v_i(t)}
\newcommand{\acceli}[0]{\dot{v}_i(t)}
\newcommand{\thrusti}[0]{\tau_i(t)}
\newcommand{\massi}[0]{m_i}
\newcommand{\alphai}[0]{\alpha_i}
\newcommand{\betai}[0]{\beta_i}
\newcommand{\mui}[0]{\mu_i}
\newcommand{\nui}[0]{\nu_i}
\newcommand{\bouyancyi}[0]{b_i}
\newcommand{\quaddragi}[0]{d_{Q_i}}
\newcommand{\lindragi}[0]{d_{L_i}}

%vectors
%\renewcommand{\bm}[1]{{\bm #1}}

%navcal vectors
\newcommand{\worldvel}[0]{ \preind{w}{}\dot{\bm{p}}}
\newcommand{\lblworldvel}[0]{ \preind{w}{}\dot{\bm{p}}_{l}}
\newcommand{\lblworldpos}[0]{ \preind{w}{}\bm{p}_{l}}

\newcommand{\beamvel}[0]{\bm{v}_{beam}(t)}
\newcommand{\dopinstvel}[0]{\preind{i}{}\dot{\bm{p}}_d(t)}
\newcommand{\dopworldvel}[0]{\preind{w}{}\dot{\bm{p}}_{d}(t)}
\newcommand{\dopworldposhat}[0]{\preind{w}{}\hat{\bm{p}}_d(t)}
\newcommand{\dopworldpos}[0]{\preind{w}{}\bm{p}_d(t)}
\newcommand{\dopworldposhatini}[0]{\preind{w}{}\hat{\bm{p}}_{d}(t_0)}
\newcommand{\dopvehvel}[0]{\preind{v}{}\dot{\bm{p}}_d(t)}

\newcommand{\cvec}[0]{\bm{c}}
\newcommand{\qvec}[0]{\bm{q}}
\newcommand{\qvecdot}[0]{\dot{\bm{q}}}

%dynamics and observer vectors
\newcommand{\xvec}[0]{\left[\begin{array}{{c}} \x \\ \xdot \end{array}\right]}
\newcommand{\xvecshort}[0]{{\bm{x}}(t)}
\newcommand{\xvecdot}[0]{\left[\begin{array}{{c}} \xdot \\ \xddot \end{array}\right]}
\newcommand{\xvecdotshort}[0]{\bm{\dot{x}}(t)}

\newcommand{\xhatvec}[0]{\left[\begin{array}{{c}} \xhat \\ \xhatdot \end{array}\right]}
\newcommand{\xhatvecshort}[0]{\bm{\hat{x}}(t)}
\newcommand{\xhatvecdot}[0]{\left[\begin{array}{{c}} \xhatdot \\ \xhatddot \end{array}\right]}
\newcommand{\xhatvecdotshort}[0]{\dot{\bm{\hat{x}}}(t)}

\newcommand{\deltaxvec}[0]{\left[\begin{array}{{c}} \Delta x(t) \\ \dot{\Delta x}(t) \end{array}\right]}
\newcommand{\deltaxvecshort}[0]{\Delta \bm{x}(t)}
\newcommand{\deltaxvecdot}[0]{\left[\begin{array}{{c}} \dot{\Delta x}(t) \\ \ddot{\Delta x}(t) \end{array}\right]}
\newcommand{\deltaxvecdotshort}[0]{\Delta \dot{\bm{x}}(t)}

\newcommand{\betavec}[0]{\left[\begin{array}{{c}} 0 \\ \beta \end{array}\right]}
\newcommand{\betavecshort}[0]{\bm{\beta}}

\newcommand{\alphavec}[0]{\left[\begin{array}{{c}} 0 \\ \alpha \end{array}\right]}
\newcommand{\alphavecshort}[0]{\bm{\alpha}}

\newcommand{\outvec}[0]{\bm{w}(t)}
\newcommand{\outvecshort}[0]{\bm{w}(t)} 
\newcommand{\outhatvecshort}[0]{\bm{\hat{w}}(t)}
\newcommand{\outtildevec}[0]{\bm{\tilde{w}}(t)}
\newcommand{\OUTtilde}[0]{\tilde{W}(t)}
\newcommand{\outdeltavecshort}[0]{\Delta \bm{w}(t)} 

\newcommand{\rvec}[0]{\left[\begin{array}{{c}} 0      \\ r_i(t)   \end{array} \right]}
\newcommand{\svec}[0]{\left[\begin{array}{{c}} s_1(t) \\ s_2(t) \end{array} \right]}
\newcommand{\nuvec}[0]{\left[\begin{array}{{c}} 0 \\ \nu \end{array}\right]}
\newcommand{\nuvecshort}[0]{\bm{\nu}}
\newcommand{\muvecshort}[0]{\bm{\mu}}
\newcommand{\muvec}[0]{\left[\begin{array}{{c}} 0 \\ \mu \end{array}\right]}

\newcommand{\measnoisevec}[0]{\bm{n}(t)}

\newcommand{\instframe}[0]{\preind{i}{}\bm{p}(t)}
\newcommand{\vehicleframe}[0]{\preind{v}{}\bm{p}(t)}
\newcommand{\localframe}[0]{\preind{l}{}\bm{p}(t)}
\newcommand{\worldframe}[0]{\preind{w}{}\bm{p}(t)}

%quadratic scalars
\newcommand{\xq}[0]{\vel |\vel|}
\newcommand{\xqshort}[0]{x_q(t)}
\newcommand{\xhatq}[0]{\dot{\hat{x}}(t) |\dot{\hat{x}}(t)|}
\newcommand{\xhatqshort}[0]{\hat{x}_q(t)}
\newcommand{\deltaxq}[0]{\dot{\hat{x}}(t) |\dot{\hat{x}}(t)| - \dot{x}(t) |\dot{x}(t)|}
\newcommand{\deltaxqshort}[0]{\Delta \dot{x}_q(t)}

%matrices
\newcommand{\Amatrix}[0]{\left[\begin{array}{{cc}} 0 & 1 \\ 0 & \mu \end{array}\right]}
\newcommand{\Ashort}[0]{\bm{A}}
\newcommand{\Cshort}[0]{\bm{C}}
\newcommand{\Pshort}[0]{\bm{P}}
\newcommand{\Lshort}[0]{\bm{L}}
\newcommand{\Qshort}[0]{\bm{Q}}
\newcommand{\eye}[0]{ I}
\newcommand{\momega}[0]{{\bm{m}}(i\omega)}
\newcommand{\ms}[0]{{\bm{m}}(s)}
\newcommand{\inmap}[0]{{\bm{b}}}
\newcommand{\outmap}[0]{\bm{c}}
\newcommand{\Rshort}[0]{\bm{R}}
\newcommand{\insttovehicle}[0]{\preind{v}{i}R}
\newcommand{\vehicletolocal}[0]{\preind{l}{v}R(t)}
\newcommand{\localtoworld}[0]{\preind{w}{l}R(t)}
\newcommand{\vehicletoworld}[0]{\preind{v}{w}R}

\newcommand{\cvecTtrans}[0]{\cvec\hspace{0.03in}^T}
\newcommand{\cvecT}[0]{\cvec}
\newcommand{\cvectrans}[0]{\cvec\hspace{0.03in}^T}
\newcommand{\bvecTtrans}[0]{\bvec\hspace{0.03in}^T}
\newcommand{\bvectrans}[0]{\bvec\hspace{0.03in}^T}
\newcommand{\bvecT}[0]{\bvec}
\newcommand{\qvectrans}[0]{\qvec\hspace{0.03in}^T}
\newcommand{\AT}[0]{A}
\newcommand{\QT}[0]{Q}
\newcommand{\PT}[0]{P}

\newcommand{\zetavec}[0]{\bm{\zeta}}

%define the gravimeter position
\newcommand{\gravposvec}[0]{\preind{w}{}\bm{p_g}}
\newcommand{\gravposx}[0]{\preind{w}{}p_{g_x}}
\newcommand{\gravposy}[0]{\preind{w}{}p_{g_y}}
\newcommand{\gravposz}[0]{\preind{w}{}p_{g_z}}
\newcommand{\gravposveclong}[0]{\left[ \begin{array}{c} \gravposx \\  \gravposy \\ \gravposz \end{array} \right]}
\newcommand{\gravposvecdot}[0]{\preind{w}{}\bm{\dot{p}_g}}
\newcommand{\gravposxdot}[0]{\preind{w}{}\dot{p}_{g_x}}
\newcommand{\gravposydot}[0]{\preind{w}{}\dot{p}_{g_y}}
\newcommand{\gravposzdot}[0]{\preind{w}{}\dot{p}_{g_z}}
\newcommand{\gravposvecddot}[0]{\preind{w}{}\bm{\ddot{p}_g}}
\newcommand{\gravposxddot}[0]{\preind{w}{}\ddot{p}_{g_x}}
\newcommand{\gravposyddot}[0]{\preind{w}{}\ddot{p}_{g_y}}
\newcommand{\gravposzddot}[0]{\preind{w}{}\ddot{p}_{g_z}}

%define the translational offset from the vehicle frame to the gravimeter
\newcommand{\gravoffsetvec}[0]{\preind{g}{v}\bm{d}}
\newcommand{\gravoffsetx}[0]{\preind{g}{v}d_x}
\newcommand{\gravoffsety}[0]{\preind{g}{v}d_y}
\newcommand{\gravoffsetz}[0]{\preind{g}{v}d_z}
\newcommand{\gravoffsetveclong}[0]{\left[ \begin{array}{c} \gravoffsetx \\  \gravoffsety \\ \gravoffsetz \end{array} \right]}

\newcommand{\Rvehtograv}[0]{\preind{g}{v}R}
%define vehicle heading, pitch, and roll
\newcommand{\hdg}[0]{\phi}
\newcommand{\pitch}[0]{\theta}
\newcommand{\roll}[0]{\psi}
\newcommand{\hdgdot}[0]{\dot{\phi}}
\newcommand{\pitchdot}[0]{\dot{\theta}}
\newcommand{\rolldot}[0]{\dot{\psi}}
\newcommand{\hdgddot}[0]{\ddot{\phi}}
\newcommand{\pitchddot}[0]{\ddot{\theta}}
\newcommand{\rollddot}[0]{\ddot{\psi}}

%define the vehicle position
\newcommand{\vehicleposvec}[0]{\preind{w}{}\bm{p_v}}
\newcommand{\vehicleposx}[0]{\preind{w}{}p_{v_x}}
\newcommand{\vehicleposy}[0]{\preind{w}{}p_{v_y}}
\newcommand{\vehicleposz}[0]{\preind{w}{}p_{v_z}}
\newcommand{\vehicleposveclong}[0]{\left[ \begin{array}{c} \vehicleposx \\  \vehicleposy \\ \vehicleposz \end{array} \right]}
\newcommand{\vehicleposxdot}[0]{\preind{w}{}\dot{p}_{v_x}}
\newcommand{\vehicleposydot}[0]{\preind{w}{}\dot{p}_{v_y}}
\newcommand{\vehicleposzdot}[0]{\preind{w}{}\dot{p}_{v_z}}
\newcommand{\vehicleposxddot}[0]{\preind{w}{}\ddot{p}_{v_x}}
\newcommand{\vehicleposyddot}[0]{\preind{w}{}\ddot{p}_{v_y}}
\newcommand{\vehicleposzddot}[0]{\preind{w}{}\ddot{p}_{v_z}}

%define the vehicle velocity vec
\newcommand{\vehiclevelvec}[0]{\left[ \begin{array}{c} \vehicleposxdot \\  \vehicleposydot \\  \vehicleposzdot \\ \hdgdot \\ \pitchdot \\ \rolldot \end{array} \right]}
%define the vehicle accel vec
\newcommand{\vehicleaccvec}[0]{\left[ \begin{array}{c} \vehicleposxddot \\  \vehicleposyddot \\  \vehicleposzddot \\ \hdgddot \\ \pitchddot \\ \rollddot \end{array} \right]}

\section{Introduction}

Consider the nonlinear discrete-time system with potentially nonlinear process and measurement terms:
\[
x_{k+1} = g_k(x_k) + \epsilon_k, \quad z_{k} = h_k(x_k) + \mu_k \quad k = 0,1,\dots, M.
\]
Smoothing seeks to discover the set of best estimates $\{\hat x_k\}_{k=1}^M$ after making observations of the state, $\{z_k\}_{k=1}^M$. $\epsilon_k$ and $\mu_k$ denote noise. The notion of `best' informs the approach. When $g$ and $h$ are linear, and the innovation and measurement errors are assumed Gaussian; the classic Kalman smoother for linear dynamical systems is an efficient algorithm to obtain the mean and variance of the conditional Gaussian distribution of the state sequence given an observation sequence~\cite{anderson2005optimal}. 
Without these strong assumptions, the posterior distribution is intractable, and the optimization approach to smoothing instead seeks to find the maximum a posteriori (MAP) likelihood estimate~\cite{aravkin2017generalized}. Closely related to these ideas is the Extended Kalman Filter (EKF), see e.g.~\cite{ribeiro2004kalman}, which linearizes $g$ and $h$ around particular estimates and then applies linear techniques.  

The Unscented Kalman Filter/Smoother (UKF/UKS)  was developed to address criticisms of the EKF approach~\cite{julier1997new}, in effect to provide a competing alternative to the local linearization approach. The UKF/UKS approach is based on the statistical view rather than the optimization view. In a nutshell, the approach creates a set of `sigma points' whose sample mean and variance represent the distribution of the estimate $\hat x_k$, and then directly propagate these sigma points through nonlinear dynamics, using the sample means and variances of the output to represent $\hat x_{k+1}$. This creative approach is able to leave behind the `local' nature of any optimization based approaches, with different sigma points potentially capturing information from different regions of the domain. However, there is no obvious way to improve on the UKF/UKS estimate with additional computation, and to date we are not aware of methods that can boost performance of UKF/UKS from the initial estimate. In contrast, the optimization approach views the classic EKF method as just the first step of a standard process, and there is an obvious strategy to improve performance using additional computation. 

The origin of the optimization approach for nonlinear smoothing most traces to the work of~\cite{bell1993iterated}, who considered the classic EKF  
linearization technique and re-imagined it as the first iteration of a local Gauss-Newton method in obtaining an updated estimate $\hat x_{k+1}$ from a current estimate $\hat x_k$. MAP techniques for non-Gaussian models appeared earlier in the statistical literature~\cite{fahrmeir1991kalman}. The OKS approach~\cite{aravkin2013optimizationPOVKS} 
{postulates and solves} the MAP objective over the entire state space sequence:
\[
\min_x \rho_p(G(x)) + \rho_m(H(x) - z)
\]
where $G(x)$ captures all process differences, i.e. has components $[x_{k}-g_k(x_{k-1})] = G_{k}$, while $H$ captures all the measurements, i.e. has components $h_k(x_k) = H_{k}$ and  $z$ is the set of observations, $\{z_k\}$ in a vector, $z = [z_1, z_2, \cdots, z_M]^T$. The losses, $\rho_p$ and $\rho_m$, reflect assumptions on the process and measurement residuals; in this paper we stay with the Gaussian case and these are weighted least squares  
\begin{align}
    \min_{x\in \mathbb{R}^{nM}} \rho \Big( F(x) \Big) = \frac12\|G(x)\|^2_{Q^{-1}} + \frac12 \|H(x) -z\|^2_{R^{-1}} \label{eqn:Opt}
\end{align}
where $x = [x_1, \cdots, x_{M}],\,  x_i \in \mathbb{R}^n$,  $$F(x) = \begin{bmatrix}
    G(x) \\H(x)-z \end{bmatrix} $$ $$ \rho\left(\begin{bmatrix} 
    f_1\\f_2\end{bmatrix}\right) = \frac12\left( \|f_1\|_{Q^{-1}}^2 + \|f_2\|_{R^{-1}}^2\right).$$

The Gauss-Newton method for~\eqref{eqn:Opt} is subtly different from the EKF/EKS approach. While the EKF/EKS linearize each component $g_k$ and $h_k$ at the `current' estimate $\hat x_k$, the GN approach treats the entire sequence as a batch, so at each iteration we linearize $G$ and $H$ at the previous estimate of the entire state sequence 
$\hat x = [\hat x_0, \hat x_{t_1}, \cdots, \hat x_{t_M} ]$. 
In principle this could negatively impact performance, since EKF/EKS may potentially propagate information faster within these full-sequence iterations, rather than just between them, but in practice we find this does not play a role.
In fact, to the point of ~\cite{julier1997new} about the difficulty of dealing with EKS, it turned out to be far easier for us to use a single iteration of OKS as the de facto EKS baseline than to implement a third-party EKS algorithm.  

The Gauss-Newton method (GN) preserves the classic Kalman problem structure.  Taking~\eqref{eqn:Opt} and linearizing $G$ and $H$, each GN update requires a solution of a particular system:
\begin{equation}  \label{eq:tridiag}
\left(\tilde G^T Q^{-1}\tilde G + \tilde H^T R^{-1}\tilde H \right)x = \tilde G^TQ^{-1}\tilde W +\tilde H^TR^{-1}\tilde Z
\end{equation}
where $\tilde G$ and $\tilde H$ represent Jacobians of $G$ and $H$ at the current state space sequence, while $\tilde W$ and $\tilde Z$ represent the updated right hand sides. The system~\eqref{eq:tridiag} is specifically written to highlight the fact that the block tridiagonal structure of the classic linear KF/KS system is preserved~\cite{aravkin2021algorithms}. This makes it possible to solve~\eqref{eq:tridiag} efficiently with a two-pass algorithm at each iteration. 
In essence, we scale the classic computation, doing all the work required for a linear smoother each time we do a GN iteration.

We follow Algorithm 2.1 from~\cite{BurkeDescent} in defining the updates, including the suggested Armijo line search, and as such claim the same convergence guarantees, for example that the algorithm will decrease the value of the function output at every iteration.
Theorem 2.4 also provides global convergence guarantees for this method.
The radius of convergence for this method is inversely proportional to the nonlinearity of the model, thus we test with initial guesses close to the ground truth. In practice the KS process will typically take place after filtering or some other prediction method and assuming a high quality filtering process, we assume the smoothing will start close enough to converge. For rates of smoothing, reference recent works such as \cite{drusvyatskiy2017efficiencyminimizingcompositionsconvex}.

In the rest of this paper, we compare the performance of our baseline EKS (single GN iteration), the UKS, and the OKS in linear and nonlinear regimes. 
For linear regimes all of the smoothers are essentially equivalent, while for the highly nonlinear non-harmonic oscillator, the OKS has a clear advantage in the high noise setting.

\section{Experiments}
\label{sec:experiments}

In this section, we compare the Extended Kalman Smoother (EKS), Unscented Kalman Smoother (UKS) and Optimized Kalman Smoother (OKS) for linear and nonlinear dynamical systems. We first describe the models for each regime, and then present the oracle parameter search we use to make sure that each smoother is performing at its best. 

\subsection{Linear Regime}

For the first example, we use the classic model of tracking a smooth underlying signal from noisy measurements. 
We generate a sine wave, and measure it with noise as follows. \\
Using a period of $2\pi$ we generate our list of equally spaced measurement times ($\Delta t \approx 0.084s$ for $75$ points across the period) and use numpy's sine and cosine functions as the ground truth for all three state components. Measurements are direct observations of $x_t$ and $\dot x_t$ with additive Gaussian noise. 

Next we explain how we model the problem using Kalman smoothing. The process model is given by 
 \begin{align*} 
 GX_t &= X_{t+\Delta t}\\
 G = 
     \begin{bmatrix}
         1 & 0 & 0\\ 
         \Delta t & 1 & 0  \\
         0& \Delta t & 1
     \end{bmatrix} &
,      \quad X_t = \begin{bmatrix}
         \ddot x_t\\ \dot x_t\\ x_t
     \end{bmatrix}
 \end{align*}

The covariance of the process model is taken to be the variance of integrated Brownian motion across a period of $\Delta t$ scaled by an unknown factor $\sigma_p^2$, ~\cite{jazwinski2007stochastic}:

\noindent
\begin{equation} \label{eq:Qk} 
   Q_0 = \begin{bmatrix}
         1&0&0\\0&1&0\\0&0&1
     \end{bmatrix}, \,\,
   Q_k =\sigma_p^2
  \begin{bmatrix}
         \dt & \frac{\dt^2}{2} & \frac{\dt^3}{6}\\ 
         \frac{\dt^2}{2} & \frac{\dt^3}{3} & \frac{\dt^4}{8}  \\
         \frac{\dt^3}{6} & \frac{\dt^4}{8} & \frac{\dt^5}{20}
     \end{bmatrix} 
    \end{equation}

The observation model encodes the direct noisy observations of the position and velocity. This gives a linear model, $H$, shown below. We assume that measurement errors are independent, Gaussian with unknown variance $\sigma_m^2$, and have the following covariance matrix ($R_k$):

\noindent

 \begin{equation}
   H = \begin{bmatrix} 0&1&0\\0&0&1 \end{bmatrix}    \label{eq:Hmat} \end{equation} 

 \begin{equation}
    R_k= \sigma_m^2 \begin{bmatrix} 1&0\\0&1 \end{bmatrix} \label{eq:Rk} \end{equation}

\subsection{Non-Harmonic Oscillator:}

In this section, we describe the non-harmonic oscillator, a highly nonlinear 
process that evolves according to the following second order differential equation.
\begin{equation}
\label{eq:NHO}
     \frac{d^2x}{dt^2} = -\omega_0^2 x - \beta\frac{dx}{dt} + k_2x^2 +k_3x^3 +\xi(t)
\end{equation}
Here $\xi \sim \mathcal{N}(0,\sigma^2_p)$ is additive noise that impacts acceleration. In Equation \ref{eq:NHO} above, $\omega_0$ regulates the period and amplitude, $\beta$ damps growth over time, while $k_2$ and $ k_3$ regulate deviations from harmonicity/linearity, with impacts to the period, amplitude, and deviations from sinusoidal motion.
\begin{align*}
    \centering
     \begin{tabular}{ccc}
        \text{variable} & \text{value} & \text{major impact}\\ \hline
        $\omega_0$ & 5 &  \text{starting amplitude }($\sim0.1$ \text{ for } $x$)  \\
        \hline
        $\beta$ &  1.5 & \text{ damping of oscillations over time}\\
        \hline
        $k_2$ & 90 & \text{deviations from sinusoidal motion}\\
        \hline
        $k_3$& -0.5 &   {period decreases over time ($\ddot{x}$)}
    \end{tabular}
\end{align*}

 The level of impact a change in these parameters has depends on the step-size being used in output model. For the experiment below we use the following values with each major impact on the model outlined below. We use $\Delta t = 0.03$s for the discretization in the generation of the ground truth, and assume observations are taken every two steps and thus define $\Delta t = 0.06$s in the smoothing process model. 

Discretizing the ODE, we obtain the following transition model: 
\[
g(X_t, \xi(t)) = G X_t + NL(X_t, \xi(t)) 
\] % check what this looks like now removed: X_{t_\Delta t]}
\begin{align*}
   \underbrace{ \begin{bmatrix}
         1 & -\beta \dt & -\omega_0^2 \dt\\ 
         \dt & 1 & 0  \\
         0& \dt & 1
    \end{bmatrix}}_G
    \underbrace{
     \begin{bmatrix}
         \ddot{x}_t\\ \dot{x}_t\\x_t
     \end{bmatrix}}_{X_t} + \underbrace{ \begin{bmatrix}
        k_2\Delta tx_t^2 +k_3\Delta tx_t^3 +\xi(t) \\
        0 \\  0
    \end{bmatrix}}_{NL}
\end{align*}
The process model has both linear and nonlinear components, with nonlinearity coming from the link between acceleration and position. 
The noise term on acceleration propagates through the process, and we again have the process covariance model~\eqref{eq:Qk}. 

Our observations take the form $h(X_t) = [\dot{x}_t, x_t]$, 
and we assume the measurement errors are Gaussian, so we in fact have the same linear observation model~\eqref{eq:Hmat}  and covariance matrix~\eqref{eq:Rk} 
as in the linear case. 

\subsection{Oracle Parameter Search}

For each experiment, we generate a 2D grid for parameters
$\sigma^2_p$ from~\eqref{eq:Qk} and $\sigma^2_m$ from~\eqref{eq:Rk}. These grids are hand tuned to cover the ideal ranges for all 3 smoothers and to be fine enough to find ideal parametrizations.  For each parameter, we record the performance and later identify the best achievable error with respect to ground truth. 

The prior covariance matrix $Q_0$ is not scaled, so the two parameters control trade-off between both process and measurement pull, and between prior and data strength. 

Results of these parametrized smoothers and visualizations of how sensitive they are to the choice of the parameters are described in the next section.

\section{Results}\label{Sec:Res}

For each experiment described in Section~\ref{sec:experiments}, the error metric, $L^2$, is time-averaged least square errors to the true signal, also referred to as ground truth in the following section. 
\[ L^2(S) = \frac{\|S-G\|_2}{\|G\|_2} = \sqrt{\frac{\sum_i (s_i-g_i)^2}{\sum_j (g_j)^2}}\] %<- not squared

\subsection{Linear Regime}
\begin{figure}[ht]
    \centering
    \includegraphics[width=1\linewidth]{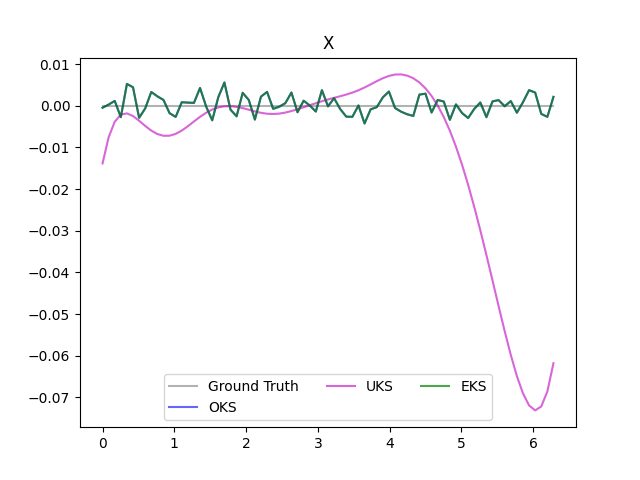}
    \caption{Deviation from Ground Truth, $\sigma_{m, \text{true}}^2 = 0.3$}
    \label{fig:LINPos}
\end{figure}

\begin{figure}[ht]
    \centering
    \includegraphics[width=1\linewidth]{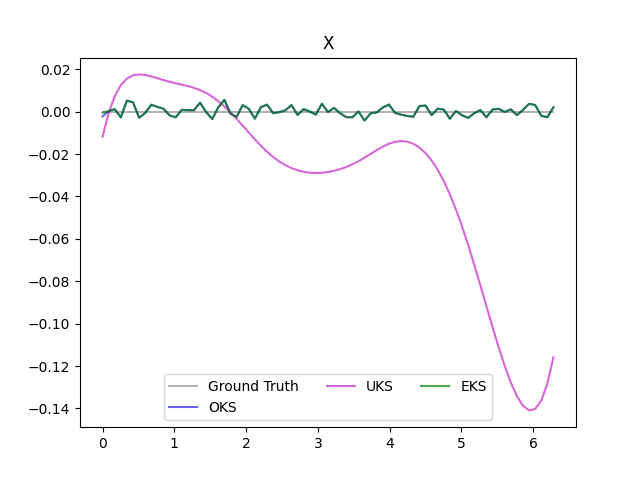}
    \caption{Deviation from Ground Truth, $\sigma_{m, \text{true}}^2 = 0.5$}
    \label{fig:LINPosNse}
\end{figure}

\begin{figure}[ht]
    \centering
    \includegraphics[width=1\linewidth]{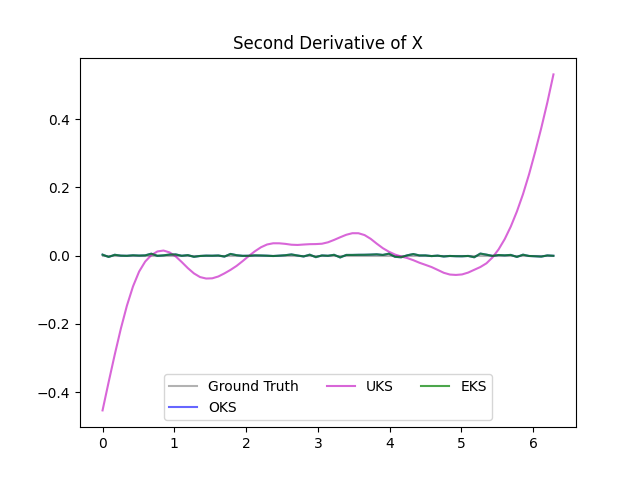}
    \caption{Deviations from Ground Truth Second Derivative, $\sigma_{m, \text{true}}^2 = 0.3$}
    \label{fig:LINSD}
\end{figure}

\begin{figure}
    \centering
    \includegraphics[width=1\linewidth]{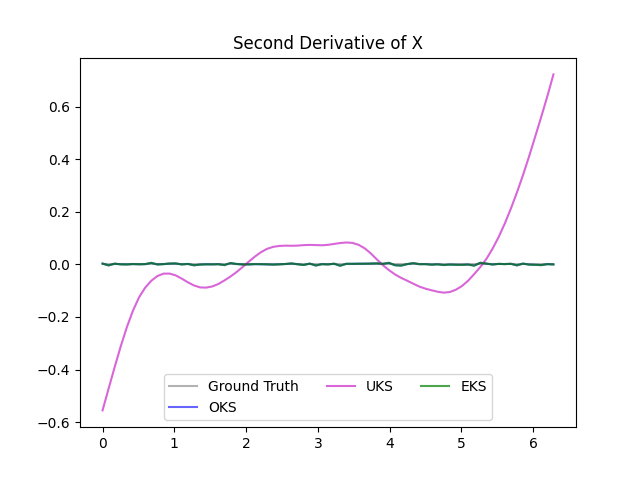}
    \caption{Deviations from Ground Truth Second Derivative, $\sigma_{m, \text{true}}^2 = 0.5$}
    \label{fig:LINSDnse}
\end{figure}

As expected, all three smoothers produce similar results and we will look closely to see the differences between them. Figures \ref{fig:LINPos} - \ref{fig:LINSDnse} show the deviations of the smoothed outputs from the ground truth for two different variances in the measurement noise. The measurements are not shown in the figure~\ref{fig:LINPos} \& \ref{fig:LINPosNse} as the scale of the deviation is smaller than the associated noise ($\sigma^2_{m,\text{true}}$).


\begin{table}[H]
    \centering
    \begin{tabular}{c||c|c|c|| p{2.5cm}}  
     ($\sigma^2_m, \sigma^2_p$) & EKS & UKS & OKS & Error Minimization \\ \hline
       (5.0, 0.2 ) & \textbf{3.259e-3} & \textbf{0.0788} & {3.259e-3} & UKS \& EKS All\\ %EXT*&UNSC all  | x
      (5.0 , 0.5 )  & 3.259e-3 & 0.1712 &\textbf{ 3.259e-3} & OKS Component\\%OPT best x  | x
      (25.0, 0.2) &   3.302e-3 & 0.0912 & 3.302e-3 & OKS Average\\%OPT best avg |x
\end{tabular}
    \caption{Relative $L^2$ Error on Sine Curve}
    \label{tab:ErrsLIN}
\end{table}
In Table~\ref{tab:ErrsLIN}, the parameters are chosen based on the minimal error found across the smoothers. All tables in this subsection contain errors from the system with larger variance in measurement noise ($\sigma_{m, \text{true}}^2 = 0.5$).  The column titled ``Error Minimization" indicates which error is optimized to obtain each set of parameters.  If the error for the parameters in question give best performance for ground truth with respect to the specific component in question, then we write ``component" in the Error Minimization column. If  the errors across all state components are minimized then we write ``average"; otherwise, ``all" indicates that both the component in question and the combined error are minimized by the parametrization. In all figures, the results shown correspond to the parameters that minimize the component error, so the bold values in the table correspond to the associated figure.

From Figure~\ref{fig:LINPos} and others in this section, we see that the extended and optimized KS perform identically, and thus the blue and green lines overlap exactly, with only slight deviations from the ground truth. Although we do find differences in the parameters the oracle search selects for these smoothers, Tables \ref{tab:ErrsLIN},~\ref{tab:ErrsLINdX}, \& \ref{tab:ErrsLINddX} show that when selecting the same parameters the errors are identical (differences in the values occurs beyond machine precision $\mathcal{O}(10^{-18})$). As previously noted we have implemented the EKS as a single iteration of the OKS implementation, and since further iterations are not necessary in the linear regime this equivalence is expected.  

The first derivative errors are similar to those from the sine curve and thus no images of the smoothed output are included; the errors for different parametrizations are presented in Table~\ref{tab:ErrsLINdX}. 
\begin{table}[ht]
    \centering
    \begin{tabular}{c||c|c|c||c}  %TBFilled IN
    ($\sigma^2_m, \sigma^2_p$) & EKS & UKS & OKS & Error Minimization \\ \hline
    (5.0, 0.2) & \textbf{3.177e-3} & \textbf{0.0865} & {3.177e-3 } & UKS \& EKS All  \\%UNSC & EXT | dx
    (5.0, 0.5)&  3.177e-3 & 0.2356 & \textbf{3.177e-3} & OKS Component \\%OPT best |dx
    (25.0,0.2) & 3.244e-3 & 0.1220 & 3.244e-3 & OKS Average \\%OPT best avg | dx
\end{tabular}
    \caption{Relative $L^2$ Error on Derivative of Sine Curve}
    \label{tab:ErrsLINdX}
\end{table}


In Figures \ref{fig:LINSD} \& \ref{fig:LINSDnse}, we can see that the smoothed result for the UKS deviates most visibly from the ground truth at the start and end of the period of observation. This deviation is larger when the variance of the noise is larger. 
Another thing of note is that the EKS and OKS deviations are significantly smaller than the UKS and unlike the UKS, the magnitude of these deviations is not affected by the change in the variance of the noise. This is due to the unscented smoother following the measurements, which can pull the curves away from the ground truth especially at the end points of our observation period. 
%Cite/reference further discussion in Conclusions/section 4?

For all three smoothers, the parametrization plays a small role in the resulting error. With the OKS and EKS there are differences of $\mathcal{O}(10^{-4})$ across the entire search space. 
 While the UKS showcases a larger variance across the search space, best seen in the UKS column of Table~\ref{tab:ErrsLINdX}, the minimal error across all state components is achieved for the same parametrization of ($5.0, 0.2$). 

 \begin{table}[H]
    \centering
    \begin{tabular}{c||c|c|c||p{2.2cm}}  
     ($\sigma^2_m, \sigma^2_p$) & EKS & UKS & OKS & Error Minimization \\ \hline
      (5.0, 0.2)& \textbf{3.688e-3} & \textbf{0.2878} & \textbf{3.688e-3} & UKS \& EKS All \newline OKS Component   \\%OPT best ddx 
      \hline
      (25.0, 0.2) & 3.688e-3 & 0.3788 & 3.3688e-3 & OKS Average \\%OPT avg | ddx
\end{tabular}
    \caption{Relative $L^2$ Error on Sine Second Derivative}
    \label{tab:ErrsLINddX}
\end{table}

 %Oscillator
\subsection{Non-Harmonic Oscillator}

\begin{figure}[ht]
    \centering
    \includegraphics[width=1\linewidth]{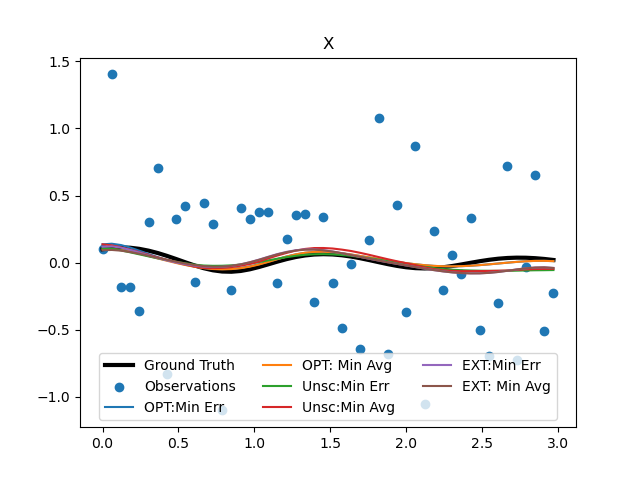}
    \caption{Oscillator Position with Noise, $\sigma^2_m = 0.5$}
    \label{fig:NHOPosNSE}
\end{figure}

\begin{figure}[ht]
    \centering
    \includegraphics[width=1\linewidth]{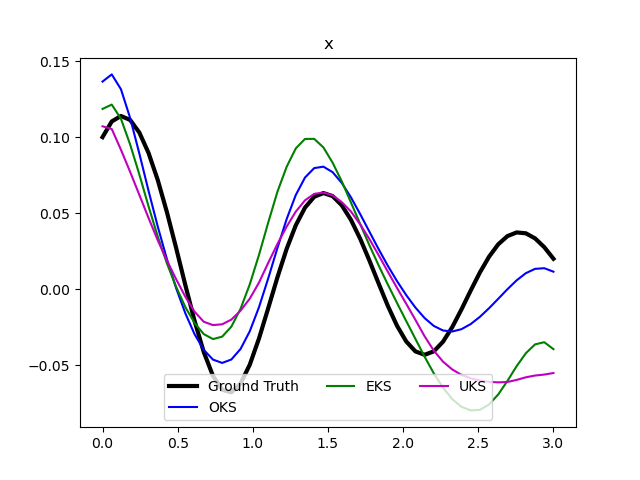}
    \caption{Oscillator Position Comparison}
    \label{fig:HOPos}
\end{figure}

\begin{figure}
    \centering
    \includegraphics[width=1\linewidth]{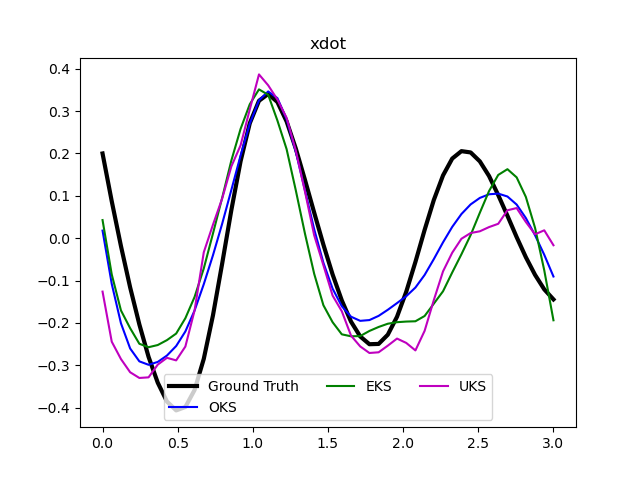}
    \caption{Oscillator Derivative Comparison}
    \label{fig:HOFD}
\end{figure}

\begin{figure}
    \centering
    \includegraphics[width=1\linewidth]{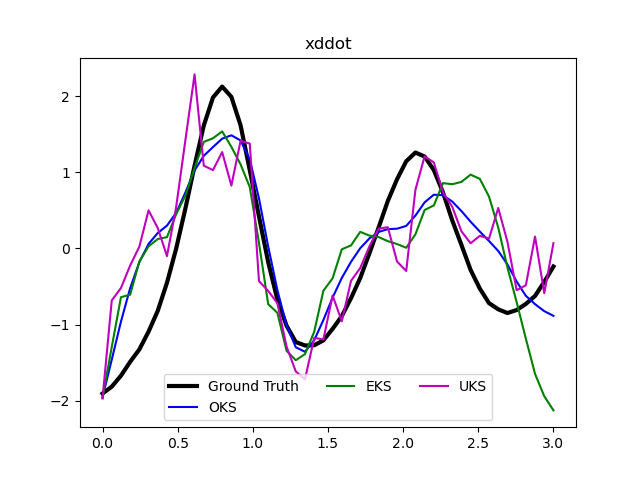}
    \caption{Oscillator Second Derivative Comparison}
    \label{fig:HOSD}
\end{figure}

In this regime, the nonlinearity of the process model causes the three smoothers to differentiate output, and gives a sense of how they behave in complex situations. In the following plots and tables, the true variance of the observations is $\sigma_{m,\text{true}}^2 = 0.5$, to highlight the benefits of the OKS in high noise regimes for nonlinear models. Figure~\ref{fig:NHOPosNSE} 
shows the observations of the position state variable compared to the ground truth and the parameter sets from Table~\ref{tab:ErrsNHO}. Figures~\ref{fig:HOPos} - \ref{fig:HOSD} show the best smoothed outputs, dropping the observations so that differences in fits vs. ground truth can be seen more clearly.


\begin{table}[ht]
\vspace{0.2cm}
    \centering
    \begin{tabular}{l||c|c|c||l}  %TBFilled IN
     ($\sigma^2_m, \sigma^2_p$)        & EKS               & UKS               & OKS               & Error Minimization \\ \hline
(0.01,1) & 0.98673 & \textbf{0.84963} & 0.50864    & UKS Min        \\
(0.025, 1) & 0.90901  & 0.90347  & \textbf{0.42640} & OKS Min        \\
(0.063, 2.51) & \textbf{0.90699} & 0.9038  & 0.4296  & EKS Min        \\
(2.51, 100)   & 0.90871 & 0.92001 & 0.44185 & OKS \& EKS Avg \\
(0.398, 100)  & 1.0883  & 0.90497   & 0.67661  & UKS Avg      
\end{tabular}
    \caption{Relative $L^2$ Error on NHO}
    \label{tab:ErrsNHO}
\end{table}

 For the nonlinear system, the smoothed output is impacted by the noise introduced at the second derivative in the process model, Equation~\eqref{eq:NHO}, when considering the position plot versus the derivatives.  The UKS, which uses statistical methods, most clearly shows this behavior. While in the position plot we see a very smooth curve, in the first and second derivative we can see the UKS line become less smooth with higher noise.
 
 In Figures~\ref{fig:HOFD} and \ref{fig:HOSD} OKS produces the smoothest output, despite the large noise. We also see that OKS, using several iterations, is able to get closer to the ground truth than either UKS or EKS. Tables~\ref{tab:ErrsNHO}, \ref{tab:ErrsNHO dx}, \& \ref{tab:ErrsNHO ddx} show the lowest achieved error, where the OKS is able to beat both the UKS and EKS. This also holds true in the case of smaller noise, but the difference between the OKS, UKS, and EKS is less pronounced.
 
\begin{table}[H]
    \centering
    \begin{tabular}{l||c|c|c||l}  %TBFilled IN
     ($\sigma^2_m, \sigma^2_p$)   & EKS   & UKS & OKS   & Error Minimization \\ \hline
(0.398, 100) & 0.83713          & \textbf{0.68821} & 0.65055          & UKS All        \\
(2.51, 100)  & \textbf{0.64782} & 0.81162           & \textbf{0.46893} & OKS \& EKS All
\end{tabular}

    \caption{Relative $L^2$ Error on NHO Dervative}
    \label{tab:ErrsNHO dx}
\end{table}

\begin{table}[ht]
\begin{tabular}{l||c|c|c||l}
($\sigma^2_m, \sigma^2_p$)  & EKS  & UKS & OKS & Error Minimization \\ \hline
(0.398, 100) & 0.91118  & \textbf{0.67025} & 0.76292 & UKS All  \\
(2.51, 100)  & \textbf{0.72228} & 0.83207 & \textbf{0.50183} & OKS \& EKS All 
\end{tabular}
\caption{Relative $L^2$ Error on NHO Second Derivative}
    \label{tab:ErrsNHO ddx}
\end{table}

\begin{figure}[ht]
    \centering
    \includegraphics[width=1\linewidth]{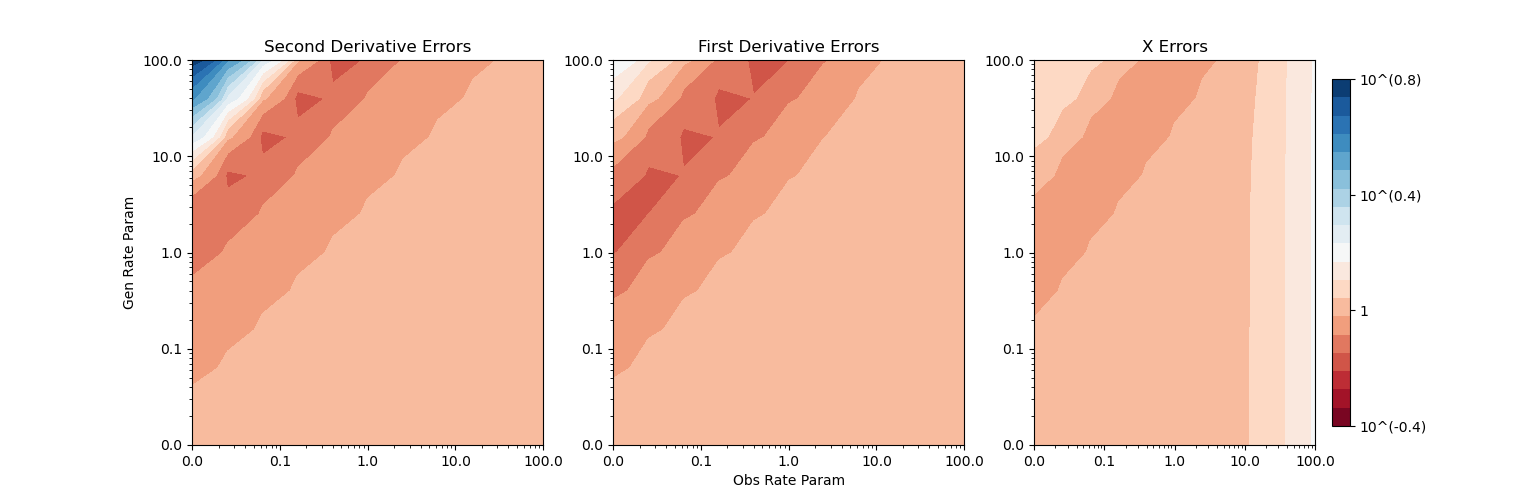}
    \caption{UKS Error Plot: Oscillator}
    \label{fig:NHO-UNSCErrPlt}
\end{figure}

\begin{figure}[ht]
    \centering
    \includegraphics[width=1\linewidth]{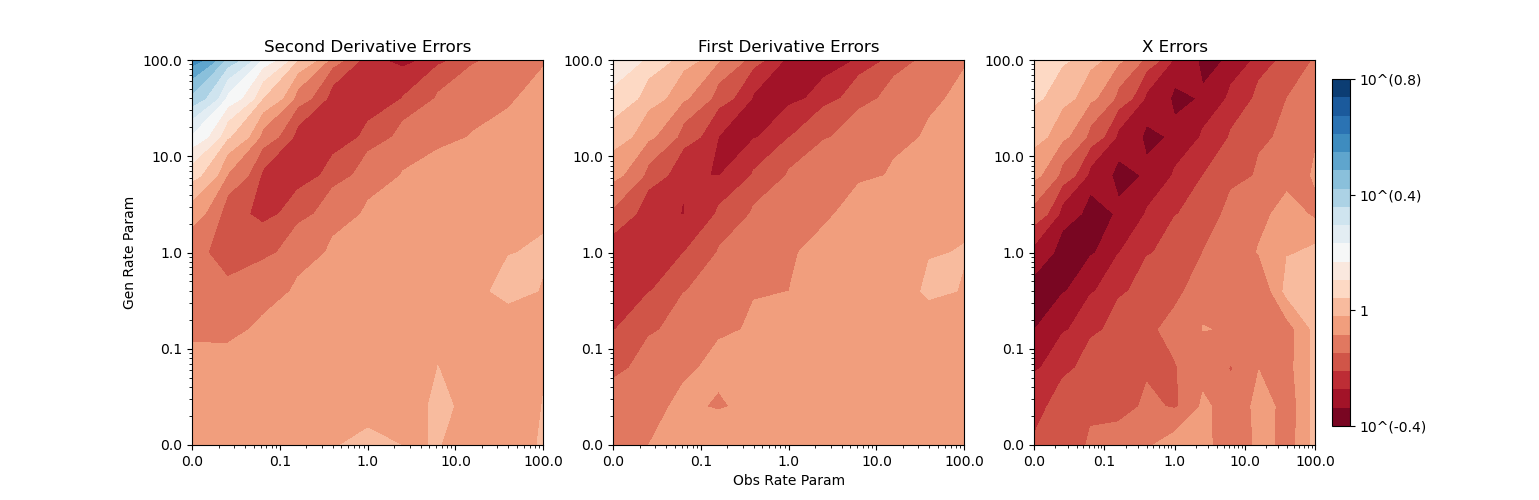}
    \caption{OKS Error Plot: Oscillator}
    \label{fig:NHO-OPTErrPlt}
\end{figure}

The sensitivity of overall performance to the two parameters in this experiment are shown in Figures \ref{fig:NHO-UNSCErrPlt} \& \ref{fig:NHO-OPTErrPlt} on a log-log scale. 
For all smoothers we can use the plots to see whether the ratio of process to measurement parameters, ($\sigma_m^2 / \sigma_p^2$), produces overlapping regions where errors across all three components is minimized. 

From Figure~\ref{fig:NHO-UNSCErrPlt}, the UKS has thin bands of minimal error for each component. These sections do not overlap significantly, which means that it is hard to find a set of parameters for UKS where all three components will be close to ground truth in the simulated example. In Table \ref{tab:ErrsNHO} above, the only parameter set where the average minimization doesn't align with the component minimization still maintains an error close to the minimal error. The average error is minimized for both the $\dot x$ and $\ddot{x}$ components, with the parameter set $(0.398, 100)$, at the top of the red band of the Figure~\ref{fig:NHO-UNSCErrPlt}.

The OKS has a larger overlapping region of minimal error.  Tables \ref{tab:ErrsNHO dx} \& \ref{tab:ErrsNHO ddx} show that the minimal component-wise error occurs under the same parameter set, $(2.51, 100)$ while for position the parameters $(0.025,1)$ have the same ratio of parameters. The plot showcases that the OKS is not sensitive to the parametrization, as the difference across the subplots is smaller than the others.\

\begin{figure}[ht]
    \centering
    \includegraphics[width=1\linewidth]{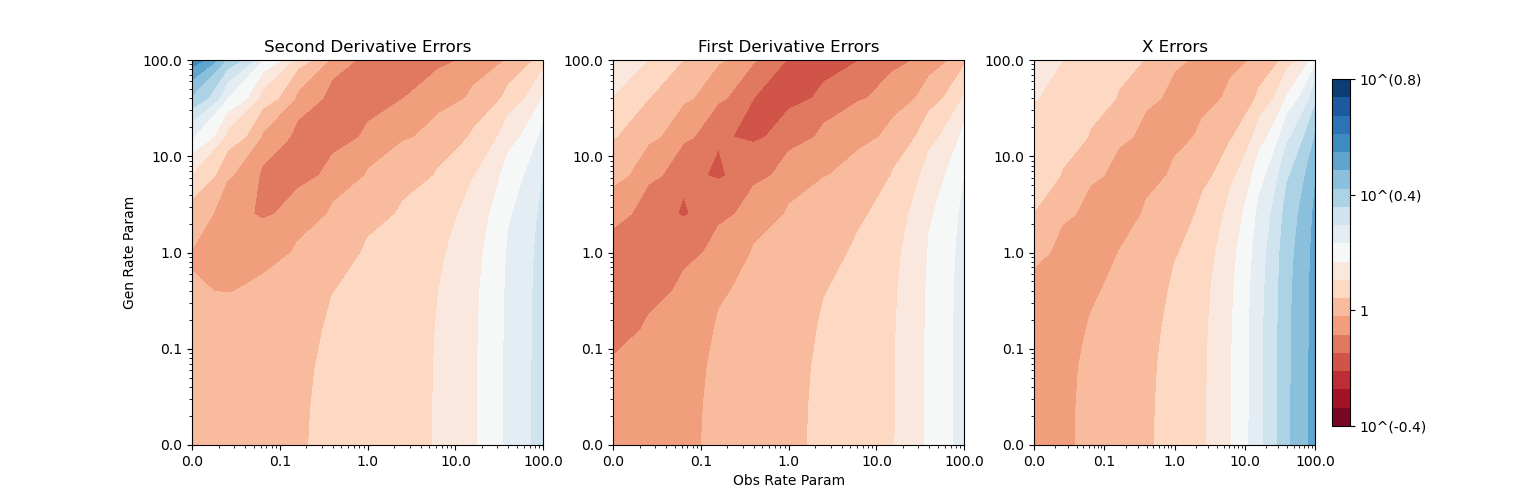}
    \caption{EKS Error Plot: Oscillator}
    \label{fig:NHO-EXTErrPlt}
\end{figure}

Figure~\ref{fig:NHO-EXTErrPlt} shows similar behavior but larger errors when the ratio of parameters is skewed. This phenomena shows EKS is more sensitive to parameters away from the optimal set, but 
does not impact the tabulated minimal error results.

\section{Conclusions}
From these results we see that although the UKS and OKS are comparable for low noise regimes (i.e. small $\sigma_m^2$), OKS is more competitive in higher error regimes. This is likely due to the difference in the optimization vs. error propagation frameworks. For example, parametrization choices scale the covariance matrices, which for the UKS directly inform the sigma points through the Cholesky decomposition of the covariance matrices done at each time step. In addition, the OKS provides an opportunity to improve performance with additional compute -- you can simply iterate the GN method -- whereas it is not clear how to do this with the UKS. 

From Figures~\ref{fig:LINSDnse} and \ref{fig:HOSD}, the UKS is more sensitive to data at the tails. The OKS and EKS also show this deviation, but appear to be more robust. 
This may be due to the ability of the OKS smoother to more efficiently pass information from the middle to the tails through the iteration process.

\section{Future Work}\label{Sec:FWork}

The comparisons in this paper use an oracle parameter search, which requires knowing ground truth in order to compare best achievable performance. In future work we plan to use data-driven techniques to compare non-oracle tuning. 
 Techniques such as  minimizing look-ahead error on a development sample can give us a more accurate understanding of the strengths and weaknesses of each method.  This is far more practical, allowing for deployment for real and complex data scenarios without need for any oracle or oracle proxy.

 The results shown are also created with the smoothers having full knowledge of the dynamics that were used to generate the ground truth. In many real systems, we do not have access to true dynamics, or  true parameters, of the systems and schemes we want to model. A related and more challenging problem is to infer ODE parameters while optimizing over the state, which would extend this work to more practical applications.

\bibliographystyle{IEEEtran}
% argument is your BibTeX string definitions and bibliography database(s)
\bibliography{swebs_refs}

\end{document}